\documentclass[%
 print,
 amsmath,amssymb,
aps,
]{revtex4-2}

\usepackage{graphicx}
\usepackage{dcolumn}
\usepackage{bm}

\renewcommand{\b}[1]{\boldsymbol{#1}}

\newcommand{\td}[2]{\frac{\mathrm{d} #1}{\mathrm{d} #2}}

\newcommand{\T}[0]{^H}

\begin{document}

\title{Data-driven linear analysis of dynamical systems via\\nonlinearity-subtracted dynamic mode decomposition}

\author{Benjamin Herrmann}%
 \email{benjamin.herrmann@uc.cl}
\affiliation{Department of Mechanical and Metallurgical Engineering \& Department of Hydraulic and Environmental Engineering, Pontificia Universidad Católica de Chile, Av. Vicuña Mackenna 4860, Santiago, Chile}

\author{Katherine Cao}
\affiliation{Department of Mechanical Engineering, Stanford University, Stanford, CA 94305, USA}

\author{Steven L. Brunton}
\affiliation{Department of Mechanical Engineering, University of Washington, Seattle, WA 98195, USA}

\author{Beverley J. McKeon}
\affiliation{Department of Mechanical Engineering \& Center for Turbulence Research, Stanford University, Stanford, CA 94305, USA}

\date{\today}

\begin{abstract}
The Dynamic Mode Decomposition (DMD) has been consolidated as a basic tool for data-driven analysis of dynamical systems, allowing simultaneous identification of coherent structures and their dynamics from time-resolved measurements.
However, with a linear regression at its core, DMD is unable to produce accurate models from recordings of dynamics that are inherently nonlinear, such as the response to large perturbations and the evolution on chaotic attractors.
Recent approaches attempt to simultaneously fit the linear and nonlinear contributions to the dynamics in a dataset by performing a regression onto a physically motivated model structure.
However, although the resulting nonlinear models can produce accurate short-term predictions, their linearization does not necessarily agree with the linearization of the original system.
In this work, we introduce a novel data-driven method — nonlinearity-subtracted DMD (NSDMD) — that focuses on producing an accurate linearization of the measured system when the nonlinear contribution to the dynamics is available while the linear part is not.
This scenario is encountered, for example, when the purely nonlinear terms in the governing equations are known, while the linear operator contains uncertain material properties or it accounts for the closure of unresolved dynamics.
Such a situation also arises when the data is generated by a black-box simulation code that is able to output the nonlinearity, but not the action of the linear operator on the snapshots.
NSDMD leverages data snapshots of the nonlinearity to explicitly account for the purely nonlinear contributions to the dynamics and formulate a regression problem that, by construction, finds a low-rank approximation of the underlying linear operator.
We demonstrate the approach on several numerical examples, showcasing its improved capabilities for data-driven linear analysis of chaotic, partially observed, advection-dominated, and high-dimensional dynamics.
\end{abstract}

\maketitle


\section{Introduction}
\label{sec:intro}

Modal decomposition techniques that rely on the inspection of a linear operator have long been a cornerstone for the analysis of spatiotemporal phenomena in structural dynamics~\citep{fu2001book}, heat transfer~\citep{bergman2011book}, and fluid mechanics~\citep{schmid2007arfm}, among many other fields in science and engineering~\citep{databook,torres2026nody}. Techniques such as linear stability, transient growth, and resolvent (input-output) analyses provide valuable insight into the dominant structures and physical mechanisms driving the dynamics of complex systems~\citep{taira2017aiaa,taira2019aiaa,herrmann2023jfm}. Moreover, modal representations have been extensively leveraged for reduced-order modeling~\citep{antoulasbook,rowley2017arfm,addison2026arxiv}.

Historically, constructing the linearized operators to be inspected has relied on explicit access to the governing equations or access to specialized, computationally intensive numerical solvers~\citep{theofilis2011arfm,frantz2023amr,kaiser2023aiaa,massaro2024cpc}. Over the past two decades, the exponential growth of high-fidelity simulations and experiments has driven a paradigm shift toward data-driven modal analysis~\citep{dmdbook,databook}. Among these approaches, the Dynamic Mode Decomposition (DMD) has emerged as a foundational tool~\citep{schmid2010jfm,rowley2009jfm}. By performing a linear regression on time-resolved state snapshots, DMD simultaneously extracts spatio-temporal coherent structures and estimates a low-rank linear operator governing their evolution, bypassing the need for intrusive access to the full discretized operators~\citep{tu2014jcd,schmid2022arfm}. Physics-informed DMD (piDMD) is a recently developed extension that allows the integration of physical principles, such as symmetries, invariances and conservation laws, into DMD models, making the approach considerably more robust to measurement noise when such priors are known~\citep{baddoo2023prsa}. Further building on DMD, \cite{herrmann2021jfm} developed an approach to perform data-driven resolvent analysis (DDRA) that is able to extract more nuanced insight into the mechanisms driving the observed dynamics. More specifically, DDRA allows the computation of resolvent modes and gains directly from data and without requiring the governing equations~\citep{herrmann2021jfm}.

An important caveat is that, because it fits a linear operator, DMD is unable to identify accurate models from data of inherently nonlinear dynamics, such as transients generated by large-amplitude perturbations or chaotic trajectories~\citep{schmid2010jfm}. As a consequence, the application of the original DDRA formulation is strictly limited to linear systems only~\citep{herrmann2021jfm}. More generally, the application of any form of data-driven linear analysis to nonlinear systems requires separating the linear and nonlinear contributions to the dynamics that drive the observed trajectories. In principle, this separation can be achieved with methods that simultaneously fit a nonlinear model that is subsequently linearized for analysis, such as the linear and nonlinear disambiguation optimization (LANDO)~\citep{baddoo2022prsa} and operator inference (OpInf)~\citep{peherstorfer2016cmame,kramer2024arfm}. However, fitting a nonlinear model with a linearization that is accurate enough for analysis may require a prohibitive amount of data if the system is high-dimensional.

\begin{figure}[t]
    \centering
    \includegraphics[width=1\linewidth]{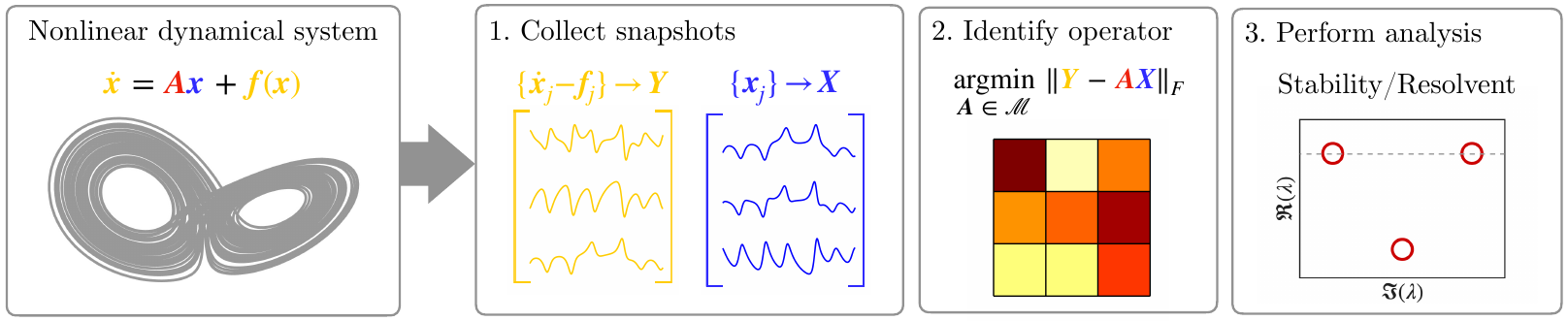}
    \caption{Schematic of the nonlinearity-subtracted dynamic mode decomposition approach (NSDMD) to perform linear analysis based on data from inherently nonlinear dynamics.}
    \label{fig:method}
\end{figure}

In this work, we introduce a novel data-driven framework---nonlinearity-subtracted dynamic mode decomposition (NSDMD)---specifically designed to produce an accurate low-rank approximation of the underlying linear operator by indirectly sampling its action on data from inherently nonlinear trajectories. By utilizing snapshots of both the state and the purely nonlinear terms in the dynamics, NSDMD explicitly subtracts out the nonlinearity to isolate the action of the linear operator, thus directly bypassing the need to fit a nonlinear model, as shown schematically in Fig.~\ref{fig:method}. We envision three distinct scenarios in which the method may be useful. First, the nonlinear contribution to the dynamics may be known a priori while the linear operator might be uncertain. This situation arises naturally in a wide class of nonlinear wave and transport equations for which the nonlinearity has a prescribed form, such as advection, while uncertain material, medium, or closure properties enter through the linear dynamics. Examples include transport processes with uncertain and perhaps spatially varying diffusion, turbulent flows with unknown linear closure terms such as eddy viscosity, and structural systems with known geometric nonlinearities but uncertain stiffness distribution. Second, the governing equations and the linear operator may in principle be known, but the available numerical infrastructure may provide easier access to the nonlinear contribution than to the action of the linear operator. For instance, an existing simulation code may provide the state fields and spatial derivatives required to evaluate an advective nonlinearity, whereas constructing or applying the corresponding linearized operator may require developing a dedicated linear solver or modifying the existing code. In this setting, NSDMD can leverage data generated by the nonlinear solver directly, while implicitly inheriting its geometry, discretization, and boundary conditions, allowing essentially the same data-driven linear-analysis procedure to be applied across different configurations. The third scenario is when
the action of the linear operator on the data snapshots can be computed directly, thus NSDMD can be applied without requiring the nonlinearity-subtraction step. In this case, the purpose of NSDMD would be to inspect the linear mechanisms that are acting on the states explored by a particular nonlinear trajectory. In the context of resolvent analysis, this enables investigating the optimal response to forcings that are present in the observed dynamics, as opposed to the most dangerous forcing, which may be unrealizable~\citep{kamal2023jfm}.

The remainder of this manuscript is organized as follows. In \S~\ref{sec:method}, we present the mathematical formulation of NSDMD along with its integration with DDRA and an extension to handle measurements on subdomains. \S~\ref{sec:examples} details the numerical examples and datasets used to demonstrate the method. \S~\ref{sec:results} presents our main results, showcasing the application of the method to chaotic, partially observed, and high-dimensional dynamical systems. Finally, in \S~\ref{sec:conclusions} we offer our conclusions and outlook on future research directions.


\section{Proposed method}
\label{sec:method}

In this section, we introduce our NSDMD method to perform data-driven linear analysis of nonlinear systems. We begin with a general description of the approach, continue with the integration of NSDMD with DDRA, and close with an extension to handle data on subdomains.

\subsection{General description}
\label{subsec:general_nsdmd}

As a starting point, we consider a generic continuous-time nonlinear dynamical system of the form
\begin{equation}
    \b{\dot{x}}=\b{Ax} + \b{f}(\b{x}),\label{dynamics}
\end{equation}
where the overdot denotes time-differentiation, $\b{x}\in\mathbb{C}^n$ is the state of the system, $\b{A}$ is the operator governing the linear part of the dynamics, and $\b{f}: \mathbb{C}^n \rightarrow \mathbb{C}^n$ corresponds to the purely nonlinear contribution to the dynamics. We point out that special attention is needed when expressing the dynamics of the system of interest in the form of eq.~\eqref{dynamics}. Importantly, we may choose to define $\b{x}$ as a perturbation from a reference state. This choice is critical since it modifies the definition of $\b{A}$ and $\b{f}$ and, therefore, it determines the interpretation of any ensuing analysis done using the linear operator. 

We now assume that a numerical representation is available for the nonlinear part of the dynamics but unavailable for the linear part. This is typically the case for black-box numerical simulations, where snapshots of $\b{f}(\b{x})$ can be outputted directly by the solver, or for low-noise and high-resolution experimental data, where it can be computed by postprocessing the state snapshots themselves. The argument for the latter is that the nonlinearity often involves local-in-space operations on the underlying discretized field and perhaps its spatial gradients, whereas $\b{A}$ is typically a denser operator that includes a Laplacian. Our aim is then to identify a low-rank approximation of $\b{A}$ from time series measurements of the state and of the nonlinearity. 

Given a set of $m$ measurements of the state $\b{x}_j=\b{x}(t_j)$ for $j= 1,\dots,m$, which may be acquired from one or several trajectories, we may assemble the data matrices
\begin{equation*}
    \b{X} = \left[\b{x}_1 \  \ \b{x}_2 \ \cdots \ \b{x}_m\right] \in\mathbb{C}^{n\times m}\quad \text{and } \quad
    \b{Y} = \left[\b{\dot{x}}_1\!-\!\b{f}_1 \ \ \b{\dot{x}}_2\!-\!\b{f}_2 \ \cdots \ \b{\dot{x}}_m\!-\!\b{f}_m\right]  \in\mathbb{C}^{n\times m},
\end{equation*}
where  $\b{f}_j=\b{f}(\b{x}_j)$ is the nonlinearity acting on $\b{x}_j$, and $\b{\dot{x}}_j=\b{\dot{x}}(t_j)$ are the time derivatives that may be approximated from the sequential state data, via finite differences, for example, if not directly measured. We point out that, in a case where we have access to the action of the linear operator on the state, we would build $\b{Y}$ directly using the snapshots of $\b{Ax}_j$, instead of indirectly via the nonlinearity-subtraction step. Once we have assembled the snapshot matrices, we may formulate an optimization problem to identify an approximation of the underlying operator $\b{A}$ that is active on the observed trajectory
\begin{equation}
  \b{A} = \underset{\b{A}\in \mathcal{M}}{\mathrm{argmin}} \ \| \b{Y}  -\b{A} {\b{X}} \|_F,
  \label{opt_nsdmd}
\end{equation}
which is identical to the piDMD problem~\citep{baddoo2023prsa}, with $\mathcal{M}$ being a matrix manifold determined by the specific DMD variant used. We remark that the choice of DMD variant used to approximate the solution to eq.~\eqref{opt_nsdmd} can be swapped as convenient. Nevertheless, for all the results presented in this paper, except where explicitly mentioned, we use the original DMD formulation~\citep{schmid2010jfm} and, therefore, $\mathcal{M}$ is the manifold of rank-$r$ matrices, for a specified value of $r$. For this case, the solution to eq.~\eqref{opt_nsdmd} is approximated by
\begin{equation}
    \b{A}=\b{Y}\b{V_{X}\Sigma}^{-1}_{\b{X}}\b{U}\T_{\b{X}},\label{dmd}
\end{equation}
where $H$ denotes the Hermitian transpose, and $\b{U_X}\in\mathbb{C}^{n\times r}$, $\b{\Sigma_X}\in\mathbb{R}^{r\times r}$ and $\b{V_X}\in\mathbb{C}^{m\times r}$ correspond to matrices containing the leading-$r$ left singular vectors, singular values and right singular vectors of $\b{X}$, respectively. We point out that the expression in eq.~\eqref{dmd} results in a large $n\times n$ matrix that is never used in practice. Instead, a much smaller $r \times r$ matrix is built by projecting onto $\b{U_X}$
\begin{equation}
    \b{\tilde{A}}=\b{U}\T_{\b{X}} \b{A}\b{U_X}=\b{U}\T_{\b{X}} \b{Y}\b{V_{X}\Sigma}^{-1}_{\b{X}}.\label{Ahat}
\end{equation}
Because $\b{\tilde{A}}$ is an $r \times r$ matrix, its eigendecomposition $\b{\tilde{A}}\b{\tilde{V}}=\b{\tilde{V}}\b{\Lambda}$ can be computed using direct solvers to obtain the matrices $\b{\Lambda}\in\mathbb{C}^{r\times r}$ and $\b{V}=\b{U_{X}}\b{\tilde{V}}\in\mathbb{C}^{n\times r}$ containing the NSDMD eigenvalues and NSDMD modes, respectively. As for the DMD problem, it can be shown that $\b{A}$ and $\b{\tilde{A}}$ share the same non-zero eigenvalues, and that the modes in $\b{V}$ correspond to eigenvectors of $\b{A}$ projected onto the column space of $\b{U_{X}}$~\citep{tu2014jcd}.

Importantly, NSDMD eigenvalues and NSDMD modes, by construction, approximate eigenvalues and eigenmodes of the underlying operator $\b{A}$. This decomposition can then be used to perform a stability analysis, or a transient growth or resolvent analysis within the DDRA framework~\citep{herrmann2021jfm}. The method is schematically summarized in Fig.~\ref{fig:method}. 

\subsection{Data-driven resolvent analysis}
\label{subsec:ddra}

Here we give a brief description of the DDRA method to comment on its integration with NSDMD. First, conventional (equation-based) resolvent analysis inspects the SVD of the operator
\begin{equation}
    \b{H}(\omega)=(i\omega \b{I} - \b{A})^{-1},\label{H}
\end{equation}
to identify the most responsive harmonic forcing, its energy gain, and the resulting harmonic response at a given frequency $\omega$ according to the linear amplification mechanisms encoded in $\b{A}$~\citep{trefethen1993science,schmid2007arfm,mckeon2010jfm}. In the absence of $\b{A}$, but equipped with an approximation for its eigendecomposition, DDRA instead inspects the SVD of the compressed resolvent operator~\citep{reddy1993jfm,herrmann2021jfm}
\begin{equation}
    \b{\tilde{H}}(\omega)=\b{F}(i\omega \b{I} - \b{\Lambda})^{-1}\b{F}^{-1} =\b{\tilde{\Psi}}(\omega)\b{\Sigma}(\omega)\b{\tilde{\Phi}}\T(\omega),\label{Hhat}
\end{equation}
where $\b{F}\in\mathbb{R}^{r\times r}$ is the Cholesky factor of $\b{V}\T\b{V}$ that is required to account for the physical inner product in eigen-coordinates, and $\b{\tilde{\Psi}}\in\mathbb{C}^{r\times r}$, $\b{\Sigma}\in\mathbb{R}^{r\times r}$ and $\b{\tilde{\Phi}}\in\mathbb{C}^{r\times r}$ correspond to the left singular vectors, singular values and right singular vectors of $\b{\tilde{H}}$, respectively. Subsequently, an approximation of the forcing and response modes that would be obtained from the full resolvent operator $\b{H}$ is given by the columns of
\begin{equation}
   \b{\Phi}(\omega)=\b{V}\b{F}^{-1}\b{\tilde{\Phi}}(\omega), \text{ and } \b{\Psi}(\omega)=\b{V}\b{F}^{-1}\b{\tilde{\Psi}}(\omega),
\label{Hhat_svd}
\end{equation}
respectively, while the associated resolvent gains are directly given by the diagonal entries in $\b{\Sigma}(\omega)$. We point out that a data-driven transient growth analysis can be carried out analogously by inspecting the compressed time-propagator operator $\b{\tilde{G}}(t)=\b{F}e^{t\b{\Lambda}}\b{F}^{-1}$ that approximates the projection of the full propagator $\b{G}(t)=e^{t\b{A}}$ onto the NSDMD modes.

It is quite clear from this pipeline that the NSDMD method seamlessly enables the use of DDRA on data from inherently nonlinear dynamics. An important subtlety is that the resulting modes and energy gains or amplifications will depend on the data snapshots used, hence providing insight into the linear amplification mechanisms that are active over the region of state space visited by the measured trajectories. Therefore, the resulting analysis might be used to answer more nuanced questions regarding the linear amplification that is present on specific nonlinear trajectories.

\subsection{Application to subdomains}
\label{subsec:subdomain}

When the system being studied, with dynamics of the form of eq.~\eqref{dynamics}, represents a spatially discretized PDE, performing an equation-based linear analysis requires careful attention to the boundary conditions, especially for advection dominated systems~\cite{cavalieri2023aiaa}. However, for a data-driven linear analysis, one would like to be able to focus on a spatial region of interest and only require measurements acquired there. For the case of numerical data, one might want to investigate the mechanisms driving the dynamics on a specific subdomain within a much larger computational domain. Moreover, experimental data of spatiotemporal phenomena is always limited to the measured field of view, that can be regarded as a subdomain with dynamics that are coupled to the unmeasured states. Therefore, a data-driven linear analysis method applied to a subdomain needs to, somehow, account for the coupling with the outside dynamics. Here we develop an NSDMD extension that is able to achieve this for the specific case of advection-dominated systems. We also explain what is the difficulty encountered for the general case, and why the proposed extension is not suitable in this scenario.

To formalize this general scenario, consider a partition of the state vector into $\b{x}=(\b{x_s},\b{x_b},\b{x_o})^{T}$, where $\b{x_s}$ denotes the state variables in the subdomain of interest, $\b{x_o}$ correspond to the unmeasured state variables outside the subdomain, and $\b{x_b}$ are state variables on a boundary region that separates the latter two in space. The dynamics of the system can then be expressed as
\begin{equation}
   \td{}{t}\begin{bmatrix} \b{x_s} \\ \b{x_b} \\ \b{x_o} \end{bmatrix} = 
   \begin{bmatrix} \b{A_{ss}} & \b{A_{sb}} & \b{0}\\
                    \b{A_{bs}} & \b{A_{bb}} & \b{A_{bo}}\\
                    \b{0} & \b{A_{ob}} & \b{A_{oo}}
    \end{bmatrix}
   \begin{bmatrix} \b{x_s} \\ \b{x_b} \\ \b{x_o} \end{bmatrix} + \b{f}(\b{x}),
    \label{eq:dynamics_partition}
\end{equation}
where the block matrices $\b{A_{ij}}$ represent the dynamics coupling the respective components of the partitioned state vector. Notice that, in eq.~\eqref{eq:dynamics_partition}, the dynamics $\b{x_s}$ and $\b{x_o}$ are uncoupled. Assuming that the underlying PDE has no non-local effects, this can always be achieved by defining $\b{x_b}$ large enough to absorb the coupling between the subdomain and the outside.

Concentrating on the subdomain only, the dynamics of $\b{x_s}$ are given by
\begin{equation}
    \b{\dot{x}_s}= \b{A_{ss}}\b{x_s} + \b{A_{sb}}\b{x_b} + \b{f_s} = \b{\Gamma}
    \begin{bmatrix}
        \b{x_s}\\ \b{x_b}    
    \end{bmatrix}
    + \b{f_s},\label{eq:dynamics_s}
\end{equation}
where $\b{\Gamma}=\begin{bmatrix}\b{A_{ss}} \ \b{A_{sb}}\end{bmatrix}$ and $\b{f_s}=\b{f}(\b{x})|_{\b{s}}$ denotes the part of the nonlinearity that contributes to the dynamics on the subdomain. From eq.~\eqref{eq:dynamics_s}, we can see that the boundary states act as exogenous inputs within the subdomain. Therefore, we can formulate a data-driven regression problem to find a low-rank approximation for $\b{\Gamma}$ by combining NSDMD with the principles used by the DMD with control (DMDc) method~\citep{proctor2016jads}.

First, we define our data snapshots as
\begin{equation*}
    \b{X} = \begin{bmatrix}\b{x}_{\b{s}1} \  \ \b{x}_{\b{s}2} \ \cdots \ \b{x}_{\b{s}m}\\
    \b{x}_{\b{b}1} \  \ \b{x}_{\b{b}2} \ \cdots \ \b{x}_{\b{b}m}
    \end{bmatrix}, \quad \text{and} \quad 
    \b{Y} = \begin{bmatrix}\b{\dot{x}}_{\b{s}1}\!-\!\b{f}_{\b{s}1} \ \ \b{\dot{x}}_{\b{s}2}\!-\!\b{f}_{\b{s}2} \ \cdots \ \b{\dot{x}}_{\b{s}m}\!-\!\b{f}_{\b{s}m}
    \end{bmatrix},
\end{equation*}
where $\b{X}$ is aggregating measurements on the subdomain and on the boundary region, whereas $\b{Y}$ contains the nonlinearity-subtracted rate of change in the subdomain only. An important detail is that, in order to compute the snapshots of $\b{f_s}$ one might need some entries in the snapshots of $\b{x_b}$. This will be the case if, for example, the nonlinearity involves a spatial gradient that needs to be numerically approximated. Nonetheless, notice that we are not using the outside state variables $\b{x_o}$.

The optimization problem to identify the underlying operator $\b{\Gamma}$ then reads
\begin{equation}
  \b{\Gamma} = \underset{\mathrm{rank}(\b{\Gamma})\le r}{\mathrm{argmin}} \ \| \b{Y}  -\b{\Gamma} {\b{X}} \|_F,
  \label{opt_nsdmdc}
\end{equation}
and an approximation to its optimal solution is given by
\begin{equation}
    \b{\Gamma}=\begin{bmatrix}\b{A_{ss}} \  \b{A_{sb}}\end{bmatrix}=\b{Y}\b{V_{X}\Sigma}^{-1}_{\b{X}}\b{U}\T_{\b{X}} = 
    \b{Y}\b{V_{X}\Sigma}^{-1}_{\b{X}}\begin{bmatrix}\b{U}\T_{\b{Xs}} \ \b{U}\T_{\b{Xb}}\end{bmatrix},\label{nsdmdc}
\end{equation}
where, again, we use the rank-r truncated SVD of $\b{X}=\b{U_X\Sigma_X V}\T_{\b{X}}$, and use $\b{U_{Xs}}$ and $\b{U_{Xb}}$ to denote the rows of the left singular vectors corresponding to the subdomain and boundary states. The projected operator governing the internal dynamics within the subdomain is given by
\begin{equation}
    \b{\tilde{A}_{ss}}=\b{U}\T_{\b{Xs}} \b{A_{ss}}\b{U_{Xs}}=\b{U}\T_{\b{Xs}} \b{Y}\b{V_{X}\Sigma}^{-1}_{\b{X}},\label{Asshat}
\end{equation}
and, finally, from its eigendecomposition $\b{\tilde{A}}_{ss}\b{\tilde{V}}=\b{\tilde{V}}\b{\Lambda}$ we obtain the NSDMDc eigenvalues $\b{\Lambda}$ and can compute the NSDMDc modes $\b{V}=\b{U_{Xs}}\b{\tilde{V}}$. This is slightly different to the DMDc solution that uses a projection onto the column space of $\b{Y}$~\citep{proctor2016jads}, whereas, here, we choose to project onto the subspace defined by the $\b{X}$ snapshots.

At this stage, it is important to pause and think about what possible insight about the dynamics of the system one might hope to obtain using only data in the subdomain and boundary regions. A reasonable initial expectation might be to understand the linear dynamics of perturbations supported within the subdomain. To explore this possibility, we will inspect the form taken by the resolvent operator $\b{H_{ss}}(\omega)$ that maps, in the frequency domain, exogenous disturbances supported within the subdomain to responses also supported within the subdomain. Taking the linear dynamics in eq.~\eqref{eq:dynamics_partition}, including an exogenous input only acting on the $\b{x_s}$ degrees of freedom, taking the Fourier transform, and after some algebraic manipulations, one arrives at the expression
\begin{equation}
    \b{H_{ss}}(\omega) = 
    \left\lbrace
        i\omega\b{I}-\b{A_{ss}} - \b{A_{sb}}
        \left[
            i\omega\b{I}-\b{A_{bb}} -\b{A_{bo}}
            \left(
                i\omega\b{I}-\b{A_{oo}}
            \right)^{-1}\b{A_{ob}}
        \right]^{-1}\b{A_{bs}}
    \right\rbrace^{-1}.\label{Hss}
\end{equation}

Unfortunately, eq.~\eqref{Hss} makes it evident that, even if disturbances are restricted to the subdomain, the coupling through the boundary may eventually drive the dynamics of the outside degrees of freedom, which can then feed back into the subdomain and affect the dynamics of the response. Therefore, even though we can approximate $\b{A_{ss}}$ and $\b{A_{sb}}$ from data, this is, in the general case, not sufficient to produce meaningful insight regarding the linear dynamics within the subdomain. Consequently, the development of a general method for data-driven linear analysis from partial measurements remains an open challenge.

Nevertheless, a very relevant special case is that of advection-dominated systems. If there is a well-defined advection direction, then we may assume that state variables will be strongly coupled to their downstream neighbors while weakly coupled to the upstream ones. The important consequence is that, if the $\b{x_b}$ states are located at the upstream boundary of the subdomain, then $\b{A_{bs}}\approx\b{0}$. Therefore, for an advection dominated system, the subdomain resolvent operator in eq.~\eqref{Hss} simplifies to
\begin{equation}
    \b{H_{ss}}(\omega)=\left( i\omega\b{I}-\b{A_{ss}} \right)^{-1},
\end{equation}
which can be approximated using the NSDMDc modes $\b{V}$ and eigenvalues $\b{\Lambda}$ obtained from eq.~\eqref{Asshat}, following the exact same DDRA procedure described in the previous subsection. This means that, for advection-dominated systems, NSDMDc in conjunction with DDRA allow using data on a subdomain (and an upstream boundary region) to investigate linear amplification mechanisms among forcings and responses restricted to that subdomain.


\section{Numerical examples and datasets}
\label{sec:examples}

\begin{figure}[t]
    \centering
    \includegraphics[width=1\linewidth]{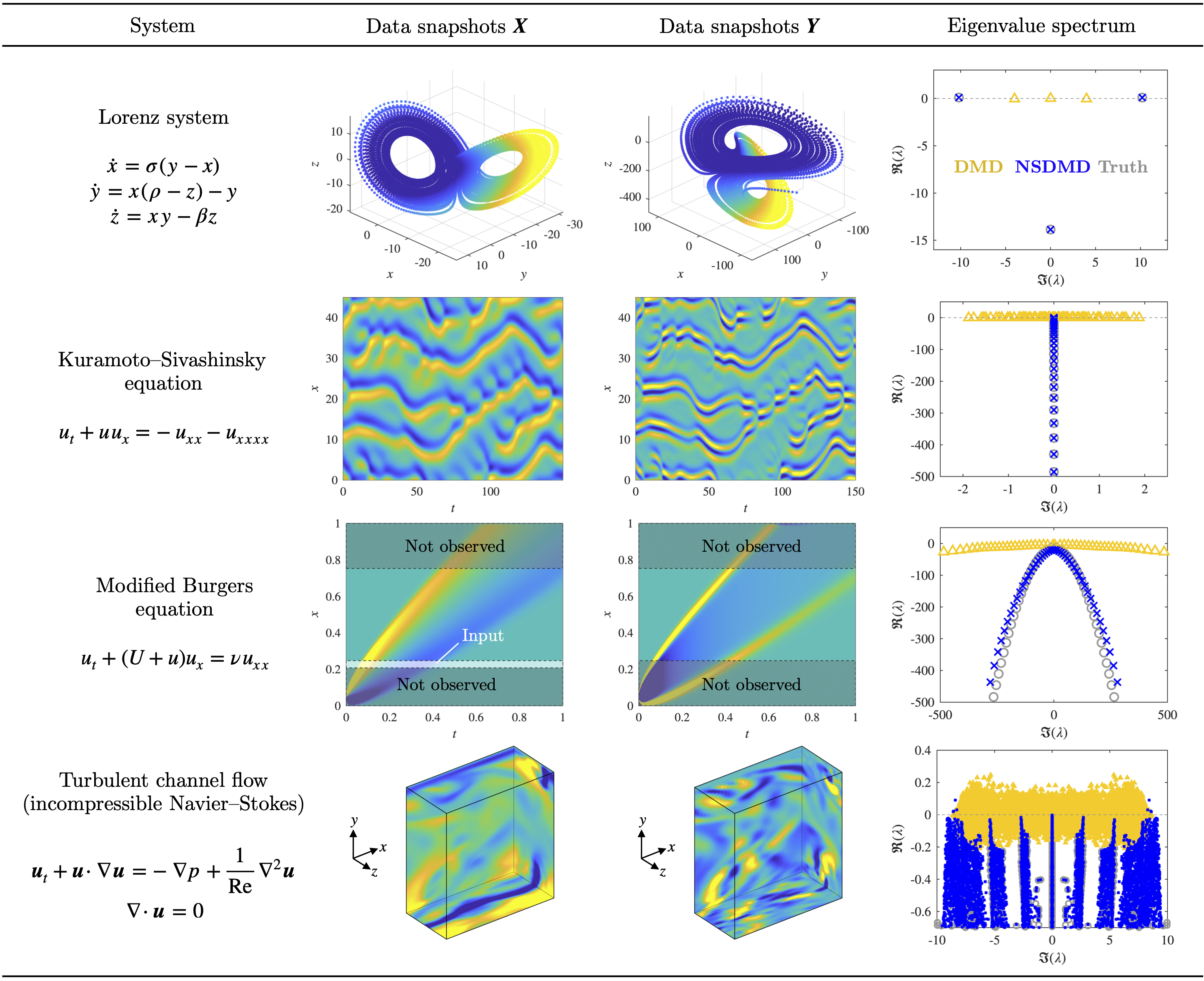}
    \caption{Overview of the application of NSDMD to four nonlinear dynamical systems. The governing equations and sample data snapshots are shown for each system, along with a comparison between the DMD, NSDMD, and true eigenvalue spectra. For the Burgers case, DMDc and NSDMDc considering boundary forcing on a subdomain are used instead. The snapshots are colored according to the fluctuations of the state variables, referenced to a base state, and their nonlinearity-subtracted rate of change. The vector norm and the streamwise component of these quantities are shown for the Lorenz system and the turbulent channel flow, where the base states are the left-lobe unstable equilibrium and the mean velocity field, respectively.}
    \label{fig:results}
\end{figure}

To demonstrate the application of NSDMD, we generate data using numerical simulations of four nonlinear dynamical systems: the Lorenz system, the Kuramoto--Sivashinsky equation, a modified Burgers equation, and the turbulent flow in a plane channel governed by the incompressible Navier--Stokes equations. These are selected to feature a pedagogic example, and chaotic, advection-dominated, and high-dimensional dynamics, respectively. For every system, we integrate trajectories starting from specified initial conditions and collect sequences of snapshots of the state and of the nonlinearity-subtracted rate of change of the state that are then stored in the data matrices $\b{X}$ and $\b{Y}$. Sample data snapshots for each system are shown in Fig.~\ref{fig:results}. Below we describe the numerical setups and generated datasets for the four examples in more detail.

\subsection{Lorenz system}
\label{subsec:ex_lorenz}

We first consider the classic three-dimensional Lorenz system using the parameters $\sigma=10$, $\beta=8/3$ and $\rho=28$, which exhibits chaotic dynamics. For these parameter values, the system has an unstable fixed point located at $(\bar{x},\bar{y},\bar{z})=(\eta,\eta,\rho-1)$, where $\eta=-\sqrt{\beta(\rho-1)}$, that we use as a reference state to rewrite the dynamics in the form of eq.~\eqref{dynamics}. Redefining the variables $x$, $y$ and $z$ to denote fluctuations away from this reference state, the governing equations are given by
\begin{equation}
\td{ }{t}
\begin{bmatrix}
x\\y\\z
\end{bmatrix} = 
\begin{bmatrix}
-\sigma & \sigma & 0\\
1 & -1 & -\eta\\
\eta & \eta & -\beta
\end{bmatrix}
\begin{bmatrix}
x\\y\\z
\end{bmatrix}
+
\begin{bmatrix}
0\\-xz\\xy
\end{bmatrix}
\end{equation}

To collect the snapshot data, the system is evolved starting from the initial condition $\b{x_0}=(-8, 8, 27)-(\bar{x},\bar{y},\bar{z})$ using a standard Runge–Kutta integrator in Matlab. Snapshots of the state are saved every $\Delta t = 0.005$ time units over a horizon of $50$ time units, resulting in a total of $m=10001$ snapshots that are then assembled into the $\b{X}$ data matrix. Subsequently, these are postprocessed to approximate the time derivative of the state using an eighth-order central finite difference scheme, and to compute the snapshots of the nonlinearity that are then used for assemble $\b{Y}$.

\subsection{Kuramoto--Sivashinsky equation}
\label{subsec:ex_ks}

The Kuramoto--Sivashinsky equation is a canonical one-dimensional PDE that models laminar-turbulent transition, flame propagation, and phase turbulence in reaction-diffusion systems. The equation is given by

\begin{equation}
    u_t + u u_x =- u_{xx} - u_{xxxx},
    \label{eq:ks_pde}
\end{equation}
where we have adopted the subscript notation to denote partial derivatives. We consider periodic boundary conditions and a domain size $L=45$, for which the system exhibits chaotic dynamics. A Fourier spectral collocation method is used to discretize in space using $n=512$ points, and a stiff ODE integrator implemented in Julia is used to advance the solution in time. The system is evolved starting from an arbitrary initial condition over a window $1000$ time units to reach the chaotic attractor, and then for an additional $15000$ time units over which the data snapshots are collected. The state, its time derivative, and the nonlinearity are outputted directly from the solver and saved every $\Delta t=1$ time units, leading to a total of $m=15001$ snapshots in each data matrix. Out of this long dataset, the first $151$ snapshots, over the first $150$ time units, are extracted to form a second, much shorter, dataset. Both datasets are used to compare the data requirements of NSDMD against those of other approaches in \S~\ref{subsec:res_ks}.

\subsection{Modified Burgers equation}
\label{subsec:ex_burgers}

For an example of an advection-dominated system, we consider the one-dimensional Burgers equation modified by including an additional linear advection term. The system is governed by
\begin{equation}
    u_t + (U+u) u_x = \nu u_{xx},
    \label{eq:burgers_pde}
\end{equation}
where $\nu = 0.005$ is the diffusion coefficient, and $U = 1$ is the added constant advection velocity. Equation~\eqref{eq:burgers_pde} is equivalent to the standard Burgers equation, but with the field $u(t,x)$ redefined to represent perturbations of the uniform solution $\bar{u}(t,x)=U$, which is a stable equilibrium of the system for the case with homogeneous Neumann boundary conditions. We consider the domain $x\in[0,1]$, with homogeneous Dirichlet boundary conditions for eq.~\eqref{eq:burgers_pde}. The system is discretized into $n=256$ spatial points using a spectral collocation method based on a discrete sine transform and integrated in time using a standard Runge–Kutta integrator in Matlab.

A single trajectory is simulated starting from a localized wave initial condition, given by
\begin{equation}
    u(0,x)=A\sin\left(\frac{2\pi(x-x_0)}{4\sigma}\right)\exp\left(\frac{-(x-x_0)^2}{2\sigma^2} \right),
\end{equation}
that is, a sine wave modulated by a Gaussian, where $A=10$, $x_0=0.05$, and $\sigma=0.015$ are the amplitude, center location, and standard deviation for the Gaussian respectively, and the sine function has a wavelength of $4\sigma$. The system is evolved over a window of $1$ time units, and snapshots of the state, its time derivative, and the nonlinearity are outputted directly from the solver and saved every $\Delta t=0.001$ time units, leading to a total of $m=1001$ snapshots in each data matrix.

To assess the performance of the NSDMDc formulation for data on subdomains described in \S~\ref{subsec:subdomain},
we define a subdomain of the full computational domain that is given by $x\in[0.25,0.75]$. Moreover, we define an upstream boundary region given by $x\in[0.21,0.25)$ where the state data is interpreted as an exogenous input in the analysis. The part of the snapshots on the exterior of both of these regions is assumed to be unavailable and is thus discarded, as is depicted in Fig.~\ref{fig:results}.

\subsection{Turbulent channel flow}
\label{subsec:ex_turbulence}

For our last example, we consider a turbulent fluid flow governed by the incompressible Navier--Stokes equations. As discussed in \S~\ref{sec:method}, the framing of the NSDMD problem requires expressing the system of interest in the form of eq.~\eqref{dynamics}. For an incompressible turbulent flow, this can be achieved by using the mean flow $\b{\bar{u}}$ as a base state and rewriting the Navier--Stokes equations in terms of the velocity fluctuations $\b{u}'=\b{u}-\b{\bar{u}}$, leading to
\begin{align}
    \frac{\partial\b{u}'}{\partial t}= -\nabla p' -\b{\bar{u}}\cdot\nabla\b{u}'-\b{u}'\cdot\nabla\b{\bar{u}} +\frac{1}{\mathrm{Re}}\nabla^2\b{u}' - \b{u}'\cdot\nabla\b{u}' + \overline{\b{u}'\cdot\nabla\b{u}'}, \label{NS_pert}
\end{align}
where the last term contains the divergence of the so-called Reynolds stresses that arises from averaging the original equations. If the system has periodic and/or homogeneous boundary conditions for $\b{u}'$, then projecting onto a divergence-free basis makes the pressure term vanish, thus, the projected dynamics are of the form of eq.~\eqref{dynamics}, with the second, third and fourth terms on the right-hand side of eq.~\eqref{NS_pert} corresponding to the linear part, and the fifth and sixth terms comprising the nonlinear part of the dynamics. This nonlinearity is exactly what in the resolvent analysis literature~\citep{mckeon2017jfm,rolandi2024tcfd} is referred to as the nonlinear forcing, albeit here in the time domain. 

The specific configuration we consider is a turbulent channel flow following the same setup as in~\cite{herrmann2023jfm}. This corresponds to a pressure-driven turbulent flow in a doubly-periodic plane channel. The domain size is $1.83h\times 2h\times 0.92h$ length units, where $h$ is the channel half-height, along the $x$, $y$, and $z$ coordinates that indicate the streamwise (periodic), wall-normal, and spanwise (periodic) directions, respectively. For a friction Reynolds number of $Re_{\tau}=u_{\tau} h/\nu=185$, where $u_{\tau}=(\tau_w/\rho)^{1/2}$ is the friction velocity, with $\tau_w$ being the mean wall-shear stress and $\rho$ and $\nu$ the density and kinematic viscosity of the fluid, this is the smallest domain that is able to sustain turbulence and is known as a minimal flow unit \citep{jimenez1991jfm}.

We use the spectral code Channelflow \citep{gibson2008jfm,gibson2014chflow} to perform direct numerical simulations (DNS) and generate a dataset comprised of a long sequence of snapshots acquired after statistically stationarity is reached. The code uses Chebyshev and Fourier expansions of the flow field in the wall-normal and horizontal directions, and a $3^{\mathrm{rd}}$-order Adams--Bashforth backward differentiation scheme for the time integration. We find that a grid with $32\times 101\times 16$ (after de-aliasing) in $x$, $y$, and $z$ and a time step of $0.005 \, h/U$ time units, with $U$ being the centerline velocity for the laminar flow profile, are sufficient to discretize the domain and keep the CFL number below $0.55$, for the case studied. The flow is initialized from a random perturbation of the laminar profile, simulated for $440 \, h/u_{\tau}$ time units until transients have died out and statistical stationarity is reached. Then, the mean flow is computed and $75000$ snapshots of the velocity fluctuations are saved every $1.63\, \nu/u_{\tau}^2$ ($=0.2 \, h/U$) time units over an additional $1.5\times10^4 \, h/U$ time units, which is enough to get converged statistics. Snapshots for the time derivative of the velocity fluctuations are approximated using an eighth-order central finite difference scheme, and snapshots for the nonlinearity are computed using Chebyshev and Fourier spectral differentiation. After the time derivatives are computed, we down-sample the data in time, retaining $1$ of every $20$ snapshots, yielding a final total of $m=3750$ snapshots.

Subsequently, the snapshots are Fourier-transformed in the horizontal directions, so that NSDMD can be applied independently for each pair of streamwise and spanwise wavenumbers. This ensures that the resulting model respects the shift-equivariance of the flow in the homogeneous directions and in eq.~\eqref{opt_nsdmd} amounts to restricting $\b{A}$ to a manifold of block-circulant matrices~\citep{baddoo2023prsa}. Moreover, flow through a plane channel is also equivariant under the dihedral group of transformations $D_2$, meaning that, a vertical reflection, a spanwise reflection, or a rotation about the $x$-axis of an observed flow field yields another admissible flow field~\citep{sirovich1987qam}. To respect these symmetries, the appropriate transformations are applied to the Fourier-transformed fields and appended as additional data snapshots, multiplying the total number of snapshots by four, before proceeding with the NSDMD method. 

Furthermore, to compare against the results based on NSDMD, we perform resolvent analysis based on the governing equations linearized about the mean flow. The mean flow is computed from the DNS snapshots and used to build the mean-flow-linearized operator with an in-house code based on the Orr-Sommerfeld/Squire formulation. Our code uses Chebyshev spectral collocation to discretize the wall-normal direction with the same grid used in the DNS. 


\section{Results and discussion}
\label{sec:results}

In this section, we demonstrate the application of NSDMD to the four datasets described in the previous section. For every system, we directly compute the eigenvalues of the numerical representation of the corresponding linear operator $\b{A}$ that we use as ground truth. The equation-based eigenvalue spectra are compared to the data-driven approximations obtained using DMD and NSDMD, as shown in Fig.~\ref{fig:results}. Moreover, for the modified Burgers example, results shown are obtained using DMDc and NSDMDc are used instead considering the data on a subdomain with boundary inputs. To select the rank truncation of the SVDs of all data matrices $\b{X}$ in our computations, we use the optimal hard threshold criteria for singular values by~\cite{gavish2014ieeetit}.

As expected, the spectra obtained with DMD for all the examples (DMDc for the Burgers case) contains only near neutral eigenvalues that attempt to fit the inherently nonlinear trajectories and hold no relation to the underlying linearized dynamics. Moreover, eigenvalues obtained with NSDMD (NSDMDc for the Burgers case) accurately approximate the equation-based results, as shown in Fig.~\ref{fig:results}. The rest of this section discusses the results obtained to demonstrate various aspects of the proposed method. 

\subsection{Subtracting vs fitting the nonlinearity}
\label{subsec:res_ks}

\begin{figure}[t]
    \centering
    \includegraphics[width=1\linewidth]{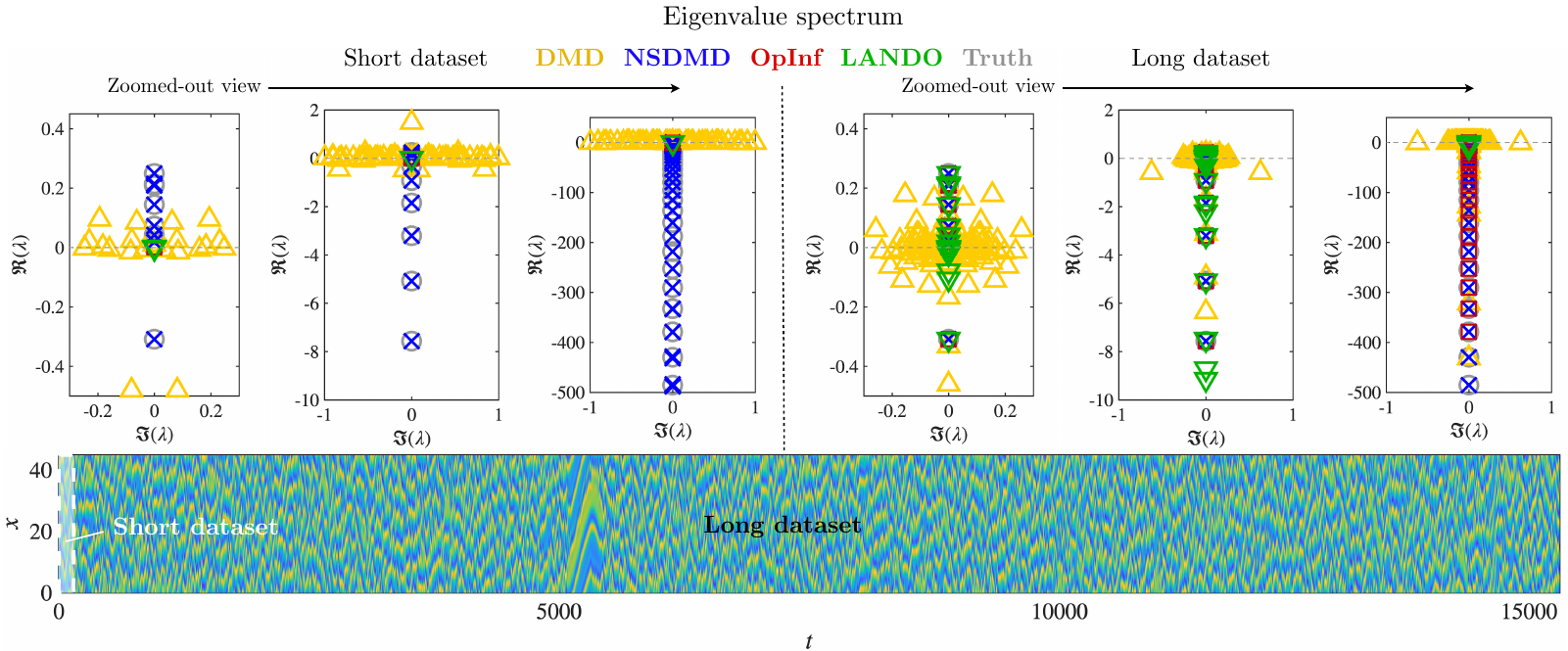}
    \caption{Comparison between the performance of DMD, NSDMD, OpInf, and LANDO to approximate the eigenvalue spectrum of the linearized Kuramoto--Sivashinsky equation from data. The methods are applied to a short dataset and a long dataset, comprised of $100$ times more data snapshots, that are also shown.}
    \label{fig:ks_results}
\end{figure}

We evaluate the performance of NSDMD on the Kuramoto--Sivashinsky equation, contrasting it directly against standard DMD and against other data-driven methods that attempt to simultaneously fit the linear and nonlinear contributions to the dynamics and later linearize the resulting model. Specifically, we use the LANDO method that uses a kernel regression on a sparse dictionary of samples~\citep{baddoo2022prsa}, and OpInf that performs an explicit regression for the POD coefficient dynamics onto a polynomial model structure~\citep{peherstorfer2016cmame}. For the LANDO implementation, we use a quadratic kernel without a bias term and do a parameter sweep to choose the dictionary sparsification tolerance $\nu$. We select the best performing case, corresponding to $\nu=1$, to be compared against NSDMD. For the OpInf implementation, we use the optimal hard threshold criteria~\citep{gavish2014ieeetit} to truncate the POD expansion and use a quadratic model structure without a bias term. The four methods, DMD, NSDMD, OpInf, and LANDO, are applied to the short and long datasets described in \S~\ref{subsec:ex_ks}, and shown in Fig.~\ref{fig:ks_results}, and compared also against the true eigenvalues of the system.

Eigenvalues obtained with all methods for both datasets are displayed in Fig.~\ref{fig:ks_results} using three levels of zoomed-out views to show different regions of the complex plane. For the short dataset, NSDMD is the only method that is able to accurately capture the true spectrum, including the unstable eigenvalues and all the stable eigenvalues deep into the negative real part of the complex plane. For the long dataset, comprised of $100$ times more data snapshots, LANDO is able to capture the unstable eigenvalues and a few of the stable eigenvalues along with a large number of spurious ones, while OpInf accurately captures most of them, and NSDMD captures all of them again. We point out that the results obtained with NSDMD applied to the short dataset are superior to those obtained with OpInf, even when the latter is applied to the long dataset. This highlights the significant impact of subtracting the nonlinearity instead of fitting it on the amount of data required to perform an accurate data-driven linear analysis. These results are explained by the difference in the amount of unknowns being identified by each regression approach, since NSDMD only requires fitting $r\times r$ coefficients for the linear dynamics, whereas OpInf needs to fit an additional $r\times r(r+1)/2$ coefficients for the quadratic terms.

\subsection{Subdomain data of advection-dominated system}
\label{subsec:res_burgers}

\begin{figure}[t]
    \centering
    \includegraphics[width=0.95\linewidth]{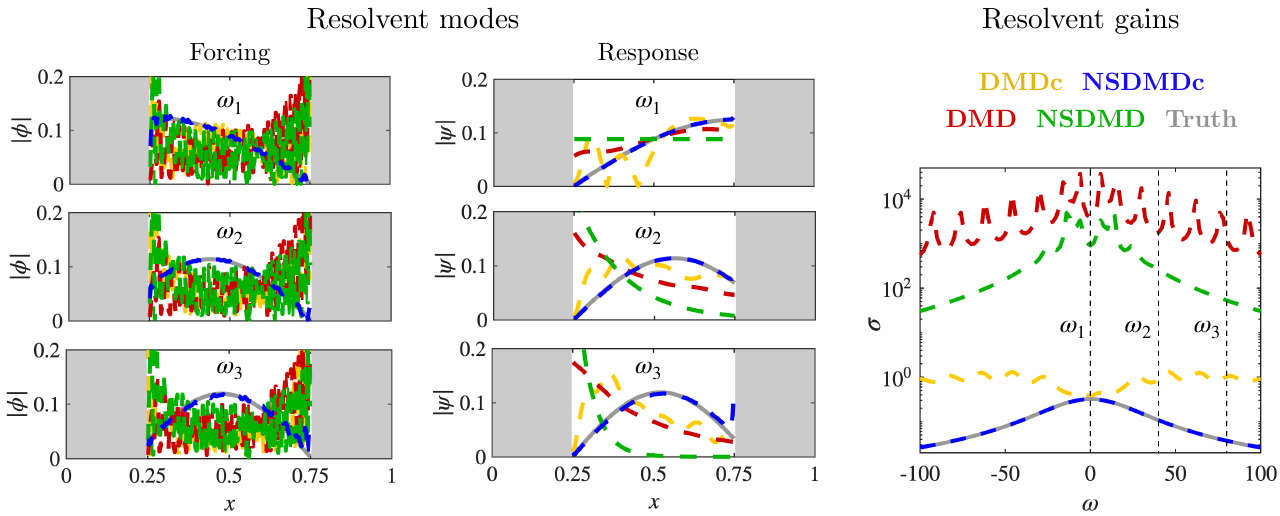}
    \caption{Application of NSDMDc in conjunction with DDRA to subdomain measurements of the modified Burgers equation. Equation-based resolvent analysis of the internal subdomain linear operator is compared to DDRA results based on the modes and eigenvalues obtained from DMD, DMDc, NSDMD, and NSDMDc. The leading resolvent gains over a range of frequencies and the leading resolvent forcing and response modes at three distinct frequencies are compared among all the approaches.}
    \label{fig:burgers_results}
\end{figure}

Here we showcase the application of NSDMDc in conjunction with DDRA on partial measurements of an advection-dominated system. We use the dataset comprised of the subdomain measurements of the modified Burgers equation described in \S~\ref{subsec:ex_burgers}. To provide a comparison against other approaches, we also apply DMD, NSDMD, and DMDc to the same dataset and then use the identified modes and eigenvalues to perform DDRA. As a reminder, here DMDc and NSDMDc treat the state data in the boundary region, highlighted in Fig.~\ref{fig:results}, as an exogenous input acting on the subdomain dynamics, whereas DMD and NSDMD use only the data on the subdomain to attempt fitting a closed-form model for the subdomain dynamics. An equation-based resolvent analysis using the subdomain linear dynamics $\b{A_{ss}}$ is performed to be used as ground truth.We compute the leading resolvent gains over a range of frequencies and the leading resolvent forcing and response modes at three distinct frequencies using the equation-based and the four data-driven approaches and compare them in Fig.~\ref{fig:burgers_results}. 

We find that, for this dataset, the only accurate DDRA results are those obtained using NSDMDc, clearly highlighting the importance of accounting for both, nonlinearity and boundary forcing. Moreover, as in~\citep{herrmann2021jfm}, we observe that forcing modes approximated from data appear noisier than response modes because $\b{U_{X}}$ is a more efficient basis to represent the latter rather than the former. We point out that this can be easily remedied with a richer dataset, whereas here we use a single trajectory to emphasize the results obtained in a data-poor regime. Finally, we remark that the presented approach is a key enabler of DDRA of nonlinear dynamics on subdomains, thus opening the door to future experimental applications.

\subsection{Application to turbulent flow}
\label{subsec:res_turbulence}

\begin{figure}[t]
    \centering
    \includegraphics[width=1\linewidth]{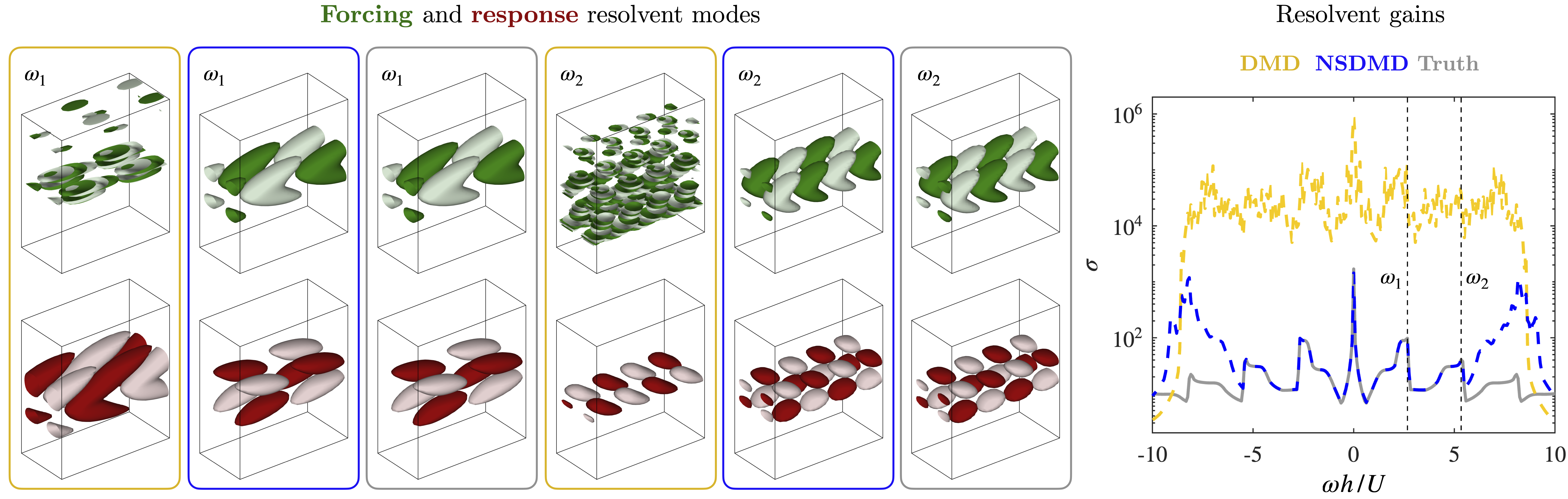}
    \caption{NSDMD-based and DMD-based DDRA of turbulent channel flow compared to equation-based resolvent analysis. Leading resolvent gain distributions and leading forcing and response resolvent modes for two distinct frequencies are shown. Forcing and response modes are visualized using isosurfaces of the wall-normal and streamwise components, respectively. Iso-values are at $\pm 0.4$ of the respective peak magnitudes, and real left- and right-going waves are added together.}
    \label{fig:chflow_results}
\end{figure}

Lastly, we present the first successful application of DDRA, enabled by NSDMD, to a turbulent flow. Specifically, we compare equation-based resolvent analysis against DMD- and NSDMD-based DDRA obtained from the turbulent channel flow dataset described in \S\ref{subsec:ex_turbulence}. We compute the leading resolvent gains over a range of frequencies and the leading forcing and response resolvent modes at two distinct frequencies using the three approaches. We find that results obtained with NSDMD match the equation-based ones, whereas those from DMD do not, as shown in Fig.~\ref{fig:chflow_results}. More precisely, we observe that NSDMD-based DDRA accurately captures the gain distribution over the frequency range $|\omega h/U|\lesssim 5.5$, outside which we get spurious amplification. This deviation is a natural limitation of the signal sampling and it should be possible to extend the range where accurate results are obtained using a smaller time step between snapshots. In fact, the observed cutoff corresponds to about a third of the Nyquist frequency. That is, we accurately identify the resolvent gains for frequencies with at least six samples per cycle, which is in line with the common practice rule of thumb used to estimate the range of frequencies to be trusted from a DMD output given the sampling rate~\citep{schmid2010jfm}. Furthermore, the leading forcing and response modes are well represented at frequencies corresponding to two peaks in the gain spectrum where structures with different streamwise wavenumbers dominate.

As previously discussed, DDRA based on a DMD model obtained from inherently nonlinear dynamics produces results that hold no connection to the underlying linear operator and their dynamical interpretation is unclear. In addition, the DMD results are highly sensitive to the amount of data and the rank truncation of the involved SVD. Therefore, it is clear that, to produce dynamically meaningful data-driven linear analysis, the nonlinearity needs to be either measured and subtracted-out or, somehow, accurately fitted to then the identified model linearized. Finally, we remark that the channel flow configuration selected is relatively simple in terms of its domain and boundary conditions, and that the application to turbulent flows over more intricate geometries or with multi-physics remains non-trivial.

\section{Conclusions}
\label{sec:conclusions}
In this work, we introduced nonlinearity-subtracted dynamic mode decomposition (NSDMD), a data-driven method designed to approximate the linearized operator underlying inherently nonlinear dynamics. By explicitly subtracting snapshots of the nonlinear terms before performing the regression, NSDMD bypasses the need to identify a complete nonlinear model and substantially reduces the amount of data required for accurate linear analysis. Across the Lorenz system, the Kuramoto--Sivashinsky equation, and turbulent channel flow, NSDMD accurately recovered the equation-based eigenvalue spectra, whereas standard DMD produces linear models that attempt to fit the observed nonlinear trajectories and fail to capture the linear mechanisms driving the dynamics. The Kuramoto--Sivashinsky example further showed that subtracting the nonlinearity can provide substantially more accurate results with considerably less data than methods that attempt to fit the linear and nonlinear dynamics simultaneously.

We also introduced NSDMDc for partial measurements of advection-dominated systems, where boundary states are treated as exogenous inputs. Combined with DDRA, this formulation accurately recovered the resolvent gains and modes of a modified Burgers equation from subdomain data, highlighting the need to account for both nonlinearity and boundary forcing. Finally, NSDMD enabled the first successful application of DDRA to a turbulent flow, producing results in close agreement with equation-based analysis. These results establish NSDMD as a practical framework for data-driven linear analysis of nonlinear systems, while extensions to noisy experimental measurements, general partially observed systems, and turbulent flows with more complex geometries remain important and challenging directions for future work.

\begin{acknowledgments}
We gratefully acknowledge P. J. Baddoo for his helpful comments and insightful discussions on the early stages of this research. This work was funded by ANID Fondecyt 1250693.
\end{acknowledgments}




\bibliographystyle{ieeetr}
\bibliography{references}

\end{document}